\documentclass[12pt]{article}
\usepackage{mathrsfs}
\usepackage{natbib,graphicx,setspace,lscape,longtable,xr}
\usepackage{mathrsfs,amsmath,amsthm,amssymb,color}
\usepackage{natbib,epsfig,graphicx,pdfpages}
\usepackage{rotating}
\usepackage{float}
\usepackage{bm}
\usepackage{CJK}
\usepackage{ragged2e}
\usepackage{setspace}
\usepackage{threeparttable}
\usepackage{multirow}
\usepackage[shortlabels]{enumitem}
\usepackage{etoolbox}
\usepackage{subfig}
\usepackage[ruled]{algorithm2e}
\usepackage{adjustbox}
\usepackage{blindtext}
\usepackage{booktabs}
\usepackage{xcolor}
\usepackage{hyperref}
\hypersetup{
	linkcolor  = blue,
	citecolor  =blue,
	urlcolor   = blue,
	colorlinks = true,
}

\bibpunct{(}{)}{;}{a}{,}{,}

\newtheorem{theorem}{Theorem}
\newtheorem{lemma}{Lemma}
\newtheorem{proposition}{Proposition}
\newcommand{\mylabel}[2]{#2\def\@currentlabel{#2}\label{#1}}
\def\beq{\begin{equation}}
	\def\eeq{\end{equation}}
\def\beqr{\begin{eqnarray}}
	\def\eeqr{\end{eqnarray}}
\def\beqrs{\begin{eqnarray*}}
	\def\eeqrs{\end{eqnarray*}}
\def\bet{\begin{theorem}}
	\def\eet{\end{theorem}}
\def\bel{\begin{lemma}}
	\def\eel{\end{lemma}}
\def\bep{\begin{proposition}}
	\def\eep{\end{proposition}}
\def\bg{\begin{figure}[tbph]\begin{center}}
		\def\eg{\end{center}\end{figure}}

\def\bc{\begin{center}}
	\def\ec{\end{center}}
\usepackage{accents}
\makeatletter
\def\widebar{\accentset{{\cc@style\underline{\mskip10mu}}}}
\def\Widebar{\accentset{{\cc@style\underline{\mskip8mu}}}}
\makeatother

\def\hat {\widehat}
\def\bar {\widebar}

\def\var{\mbox{var}}

\DeclareMathOperator*{\argmin}{arg\,min}

\def\1{\mbox{\boldmath $1$}}
\def\0{\boldsymbol{0}}

\def\mL{\mathcal L}
\def\mN{\mathcal N}
\def\mR{\mathbb R}
\def\mS{\mathcal S}
\def\mbS{\mathbb S}

\def\var{\operatorname{var}}
\def\Bias{\operatorname{Bias}}
\def\ECP{\operatorname{ECP}}
\def\MSE{\operatorname{MSE}}
\def\SE{\operatorname{SE}}
\def\CI{\operatorname{CI}}

\usepackage[a4paper,bindingoffset=0.2in,%
left=1in,right=1in,top=1in,bottom=1in,%
footskip=.25in]{geometry}

\def\boxit#1{\vbox{\hrule\hbox{\vrule\kern6pt\vbox{\kern6pt#1\kern6pt}\kern6pt\vrule}\hrule}}

\numberwithin{equation}{section}

\begin{document}
	%\begin{CJK}{GBK}{song}
	
	\begin{center}
		{\bf\Large On the asymptotic properties of a bagging estimator with a massive dataset}\\
		Yuan Gao$^{a}$,  Riquan Zhang$^{a}$\footnote{Corresponding author.\ E-mail address: rqzhang@stat.ecnu.edu.cn}\, Hansheng Wang$^{b}$\\
		{\it \footnotesize
		$^{a}$School of Statistics and KLATASDS-MOE, East China Normal University, Shanghai, China; 
		$^{b}$Guanghua School of Management, Peking University, Beijing, China 
	}
		%Yuan Gao$ ^{a} $, Xiao Wang$^{b}$, Hansheng Wang$^{c}$\\
		% {\it\small
			% 	$ ^a $ School of Statistics, East China Normal University, Shanghai, China\\
			% 	$ ^b $ Key Laboratory of Advanced Theory and Application in Statistics and Data Science - MOE, East China Normal University, Shanghai, China}

%		This version: \today
	\end{center}

\begin{abstract}
	Bagging is a useful method for large-scale statistical analysis, especially when the computing resources are very limited. 
We study here the asymptotic properties of bagging estimators for $M$-estimation problems but with massive datasets.
We theoretically prove that the resulting estimator is consistent and asymptotically normal under appropriate conditions.
The results show that the bagging estimator can achieve the optimal statistical efficiency, provided that the bagging subsample size and the number of subsamples are sufficiently large.
Moreover, we derive a variance estimator for valid asymptotic inference.
All theoretical findings are further verified by extensive simulation studies.
Finally, we apply the bagging method to the US Airline Dataset to demonstrate its practical usefulness.
\end{abstract}


	\newpage
	
	% \csection{INTRODUCTION}

\section{INTRODUCTION}

Due to the explosion of the data in this information age, massive datasets are often encountered in various fields. For example, social media tweets, online shopping records, surveillance camera videos and online search queries. They produce enormous amount of data every day \citep{fan2020statistical}.
However, the massive data size also makes the corresponding statistical estimation and inference very challenging. 
Often the data could be too large to be comfortably loaded into the computer's memory as a whole.
This means that it would be difficult, if not infeasible, to directly compute the classical whole sample based estimator. 
Consequently, effectively analyzing these massive data with statistical guarantees is a problem of great importance.

In the last decade, many methods have been proposed to deal with the large-scale statistical computation problems.
These methods can be broadly classified into two main categories.
The first is the distributed computing approach. 
The general idea behind is the divide-and-conquer, i.e, to divide a complicated problem into small pieces and then tackle them by multiple machines in a parallel way.
Distributed computing has been found very useful in many large-scale statistical problems, such as quantile regression \citep{volgushev2019distributed, chen2020distributed, pan2021note}, $M$-estimation \citep{zhang2013communication,huang2019distributed, jordan2019communication}, 
variable selection and feature screening \citep{ zhu2021least, li2020distributed}, principal component analysis \citep{fan2019distributed, chen2021distributed}, nonparametric and semiparametric regression \citep{chang2017divide,lian2019projected}.
The second approach is the subsampling method. 
The key idea is to compute the estimator based on a subsample that has been carefully selected from the whole sample. 
Consequently, the corresponding computational cost can be largely saved by choosing a relatively small sized subsample.
To take full advantage of the informative observations, different nonuniform sampling strategies have been devised, including leverage score-based subsampling \citep{drineas2012fast,ma2015statistical,ma2020asymptotic} and optimal subsampling \citep{wang2018optimal,wang2019information, wang2021optimal}. 

These subsampling methods are particularly useful when the available computing resources are very limited. 
However, two limitations of the methods are worth mentioning.
First, it requires a careful specification about the subsampling probability for each observation. This generally leads to a sampling cost of the order $O(N)$, where $N$ is the whole sample size.
Second, the resulting estimator is computed based on only one or two relatively small sized subsamples. Consequently, it is typically statistically less efficient than the whole sample based estimator.
To fix the problems, the bootstrap aggregating method \citep{breiman1996bagging, buhlmann2002analyzing}, or \textit{bagging}, provides an effective solution.
The basic idea of the bagging is to compute the estimators based on multiple subsamples generated by the method of simple random sampling with replacement. Subsequently, these estimators are aggregated into a more stable one.
Recently, similar strategies have been adopted to solve large-scale statistical problems.
For example, \cite{wu2021subsampling} proposes two subsampling-based moment estimators, whose asymptotic properties are investigated. \cite{zhu2021feature} develops a subsampling-based feature screening procedure and establishes the corresponding screening consistency.
However, the asymptotic properties of the bagging estimator seems to be largely unknown for a general $M$-estimation problem with a massive dataset. 
In this regard, we devote this work to the theoretical study of the bagging $M$-estimators with datasets of massive sizes.
We remark that our work is closely related to the existing literature but with a clear difference.
For example, \cite{wager2014asymptotic} investigates the asymptotic properties of tree-based bagging estimators (i.e., random forests). However, we focus on parametric regression models.
As another example, \cite{kleiner2014scalable} proposes the Bag of Little Bootstraps (BLB) method for distributed estimation, where weighted subsamples are used. In contrast, we use unweighted subsamples. 
This makes our method much easier to implement and study.
Recently, \cite{zou2021subbagging} studies a similar problem but with subsamples obtained by simple random sampling without replacement. However, our subsamples are obtained with replacement, which leads to a much reduced sampling cost.
This is particularly true, when the dataset is of a massive size and thus has to be placed on a hard drive. For an excellent discussion in this regard, we refer to Section 3.1 of \cite{wu2021subsampling}.

To implement the bagging method, we use the simple random sampling with replacement to generate the subsamples. 
We remark that this sampling strategy not only is simple enough to be practically implemented, but also saves the cost for specifying subsampling probabilities.
%As mentioned above, this is indeed a computationally very expensive operation.
Once a bagging subsample is obtained, we can compute the corresponding bagging subsample estimator for the parameter of interest.
After that, these subsample estimators are averaged to form the final bagging estimator.
We theoretically investigate the statistical properties of the resulting estimator and establish its consistency and asymptotic normality.
Moreover, we show that the bagging estimator can be asymptotically as efficient as the whole sample based estimator under mild conditions.
A variance estimator is also constructed for the bagging estimator to facilitate valid asymptotic inference.
Extensive simulation experiments are conducted to corroborate our theoretical findings.

The rest of the paper is organized as follows. 
Section 2 introduces the details of the problem setting and the bagging estimator. The asymptotic analysis is also included.
In Section 3, we present the simulation studies and a real data example.
The article is concluded in Section 4. 
All technical proofs are included in the Appendix.

\section{THE METHODOLOGY}

\subsection{Model and notations}

Let $Z_i  = (X_i^\top, Y_i)^\top \in \mR^{p+1}$ be the observation collected from the $i$-th subject, where $N$ is the whole sample size.  
Furthermore, $Y_i\in \mR$ is the response of interest, and $X_i=(X_{i1},\dots, X_{ip})^\top \in \mR^{p}$ is the corresponding $p$-dimensional covariate vector. 
We assume that $Z_i\ (1\le i\le N)$ are independent and identically distributed.
Let $\ell_i(\theta) = \ell(\theta; Z_i)$ be a plausible smoothing loss function. 
Then, the global loss function is defined as $\mL(\theta ) = \sum_{i=1}^N \ell_i(\theta) $, whose minimizer is given by $\hat\theta = \argmin_{\theta} \mL(\theta) $. 
Throughout this article, we refer to $\hat \theta$ as the global estimator to emphasize the fact that it is obtained based on the whole sample.
Let $\theta_0$ be the true value of the unknown parameter $\theta$, which is given by $\theta_0 = \argmin_{\theta} E\{\ell_i(\theta)\}$. In many cases,  the global estimator $\hat\theta$ admits the following asymptotic rule \citep{shao2003mathematical, van2000asymptotic},
\begin{align*}
	\sqrt{N} (\hat \theta - \theta_0) \to^d  \mN(0, \Sigma)
\end{align*}
for some positive-definite matrix $\Sigma \in\mR^{p\times p} $, as $N\to \infty$.
However, it may be difficult to compute $\hat \theta$ when the whole sample size $N$ is extremely large. 
This is particularly true if the practitioners have very limited computing resources.
Consequently, we need to have a method that should be computationally more feasible.
To this end, we find the bagging method a practically useful solution for the large-scale statistical estimation and inference.

Specifically, let $ \mbS = \{1,\dots, N\}$ denote the index set of the whole sample. 
Let $n$ be the size of the bagging subsample and $K$ be the number of subsamples. 
Write $\mS_k = \{i_1^{(k)},\dots, i_n^{(k)} \} \subset \mbS $ as the $k$-th bagging subsample, where $i_m^{(k)}$ (for any $1\le m\le n,\ 1\le k\le K$) is generated independently from $\mbS$ by the method of simple random sampling with replacement. 
In other words, conditional on $\mbS$, $i_m^{(k)}$'s are  independent and identically distributed with probability $P(i_m^{(k)} = i) = 1/N$ for any $i \in \mbS$. 
Accordingly, we can obtain an estimator based on the bagging subsample $\mS_k$ as $\hat \theta^{(k)} = \argmin_{\theta \in \mR^p} \mL_k(\theta)$ for each $1\le k \le K$, where $\mL_k(\theta) = \sum_{i\in \mS_k} \ell_i(\theta)$.
Then, we can aggregate these subsample estimators into a more stable one as $\hat \theta_\textup{BAG} = K^{-1}\sum_{k=1}^K \hat \theta^{(k)} $. 
We refer to $\hat \theta_{\textup{BAG}}$ as the bagging estimator in the following.

\subsection{Statistical properties of the bagging estimator}

To study the theoretical properties of the bagging estimator, we assume the following technical conditions. 

\begin{enumerate}[(C1)]
	\item \textbf{(Parameter Space)} The parameter space $\Theta $ is a compact and convex set in $\mR^p$. In addition, the true value $\theta_0$ lies in the interior of $\Theta$.
	
	\item \textbf{(Convexity)} The loss function $\ell(\theta; Z)$ is convex with respect to $\theta$ for almost all $Z\in \mR^{p+1}$.
	
	\item \textbf{(Gradient and Hessian)} Assume the population gradient of the loss function vanishes at $\theta_0$, that is., $E\big\{ \dot\ell_i(\theta_0) \big\} = 0$.  Let $V(\theta) =  \var\{\dot\ell_i(\theta) \}$, and write $V = V(\theta_0)$ for simplicity.
	Let $H(\theta) = E \big\{ \ddot\ell_i(\theta) \big\} $ be the population Hessian matrix of the loss function at $\theta$, and write $H=H(\theta_0)$ for simplicity.
	Assume that both $V$ and $H$ are positive-definite.
	%	Further assume that $ \lambda_{\min}(H) > \mu$ for some constant $\mu>0$, where $\lambda_{\min}(A)$ denotes the smallest eigenvalue  of an arbitrary matrix $A$.
	Furthermore, define $\Sigma = V^{-1} H V^{-1}$.
	
	\item \textbf{(Smoothness)} Define $B(\delta) = \{ \theta\in\Theta: \| \theta-\theta_0\| < \delta \}$ as a ball around the true value of $\theta_0$ with radius $\delta>0$. We next make two assumptions as follows. (C4.1) There exists some constants $G_1$ and $G_2$ such that $E\big\{ \| \dot \ell_i(\theta_0) \|^8  \big\} \le G_1^8 $ and $E\big\{ \| \ddot \ell_i(\theta_0) - H(\theta) \|^8  \big\} \le G_2^8$ for any $\theta\in B(\delta)$. (C4.2)  For almost all $Z\in\mR^{p+1} $, the Hessian matrix $\ddot \ell(\theta; Z)$ is $L(Z)$-Lipschitz continuous, i.e.,$
	\| \ddot \ell(\theta; Z) - \ddot \ell(\theta'; Z) \| \le L(Z) \| \theta - \theta'\|$ for any $ \theta,\theta' \in B(\delta)$,
	where $L(Z)$ satisfies $E \big\{  L^8(Z) \big\} \le L^8$ and $E \Big[\big\{  L(Z) - E\big( L(Z) \big) \big\}^8\Big]  \le L^8$ for some constant $L$.
	
	%	\item (Subsampling condition) The subsample size $n\to \infty$ as $N\to\infty$ and $n<N$. In addition, assume $$ 
\end{enumerate}
\noindent
Condition (C1) assumes the parameter space is compact and convex, which has been commonly used in previous studies \citep{zhang2013communication,huang2019distributed,zhu2021least}.
Conditions (C2)-(C3) are standard regularity conditions. They are commonly assumed to establish consistency and asymptotic normality for $M$-estimation \citep{van2000asymptotic, lehmann2006theory}.
Last, condition (C4) assumes that the loss function are sufficiently smooth around the true value $\theta_0$. Otherwise, the standard Taylor expansion technique cannot be applied. This is also a set of standard conditions that has been popularly used in the past literature \citep{zhang2013communication,huang2019distributed, jordan2019communication}.
We then have the following lemma.
\begin{lemma} \label{lemma: root n}
	Assume conditions (C1)--(C4). Then we have $\sqrt{n} (\hat\theta^{(k)} - \theta_0) = O_p(1)$.
\end{lemma}
%\noindent
The proof of the lemma can be found in Appendix A.
By Lemma \ref{lemma: root n}, we know that $\hat \theta^{(k)}$ is a consistent estimator for $\theta_0$ for each $1\le k \le K$. 
Furthermore, recall that $\dot\mL_k(\hat \theta^{(k)}) = 0$. 
Then by results of Lemma \ref{lemma: root n} and the mean value theorem \citep{shao2003mathematical}, we can obtain that
\begin{align*}
	0=& n^{-1} \dot\mL_k(\hat\theta^{(k)}) = n^{-1}\dot\mL_k(\theta_0) + n^{-1}\int_0^1 \ddot\mL_k\Big((1-t)\theta_0 + t \hat\theta^{(k)} \Big) dt \, (\hat\theta^{(k)} - \theta_0) \\
	=& n^{-1}\dot\mL_k(\theta_0) + H  (\hat\theta^{(k)} - \theta_0 )  +  \Delta^{(k)}, 
\end{align*}
where $\Delta^{(k)} = \Big\{ n^{-1}\ddot\mL_k(\theta_0)  - H  \Big\}  (\hat\theta^{(k)} - \theta_0 ) +n^{-1}\Big\{ \int_0^1   \ddot\mL_k\Big((1-t)\theta_0 + t \hat\theta^{(k)} \Big) dt -  \ddot\mL_k(\theta_0)  \Big\}  (\hat\theta^{(k)} - \theta_0)$.
Then we have $ \hat\theta^{(k)} = \theta_0 - H^{-1} n^{-1}\dot\mL_k(\theta_0)  - H^{-1} \Delta^{(k)}  =\theta_0 + Q_1^{(k)} + Q_2^{(k)}$.
Consequently, the bagging estimator can be represented as $\hat\theta_\textup{BAG} = \theta_0 + Q_1 +Q_2$, where $Q_1 = K^{-1} \sum_{k=1}^K  Q_1^{(k)}$ and $Q_2 =  K^{-1} \sum_{k=1}^K  Q_2^{(k)}$. 
The properties of $Q_1$ and $Q_2$ are given in the following theorem.
\begin{theorem}\label{thm:bias_var}
	Assume conditions (C1)--(C4). Then we have $\hat\theta_\textup{BAG} - \theta_0 =  Q_1 +Q_2$, where $E(Q_1) = 0$, $\var(Q_1) = \big\{ (nK)^{-1} + N^{-1} \big\} \Sigma + o\big\{ (nK)^{-1} + N^{-1} \big\}$, and $Q_2 = O_p(n^{-1})$. 
\end{theorem}

See Appendix B for the detailed proof of this theorem.
As shown in Theorem \ref{thm:bias_var}, the variability of the bagging estimator is mainly captured by the term $Q_1$, which is further determined by two parts. The first part $N^{-1}\Sigma$ is related to the whole sample size $N$, which cannot be reduced by the bagging.
The second part $(nK)^{-1}\Sigma$ is related to the bagging subsample size $n $ and the number of subsamples $K$.
This term reflects the trade-off between the sampling cost and the estimation variability.
On the other hand, the bias term $Q_2$ is mainly determined by the bagging subsample size $ n $.
This is because the averaging operation in bagging cannot reduce the estimation bias.
Therefore, in practice, we may expect that the bagging subsample size $n$ should be reasonable large (e.g., $n>\sqrt{N}$) for a satisfactory estimation accuracy.
We further establish the asymptotic normality of the bagging estimator, which is given in the following theorem.
\begin{theorem}\label{thm:anormal}
	Assume conditions (C1)--(C4) and $n / \sqrt{N} \to \infty$ as $N\to \infty$. Then we have $ \{ (nK)^{-1} + N^{-1} \}^{-1/2}  (\hat\theta_\textup{BAG} - \theta_0) \to_d \mN(0, \Sigma)$ as $N\to\infty$. 
\end{theorem}
\noindent
See Appendix C for the detailed proof of this theorem.
By Theorem \ref{thm:anormal} we can see that, as long as $n\gg \sqrt{N}$ and $nK \gg N$, the bagging estimator should be asymptotically as efficient as the global estimator.

For inference purposes, we need to estimate the variance of the bagging estimator.
Note that we have obtained $K$ estimates of $\theta_0$ based on different subsamples. 
By using these estimates, we can construct the variance estimator,
\begin{align}\label{eq:se}
	\hat\SE^2( \hat\theta_{\textup{BAG}} ) = \left( \frac{1}{nK} + \frac{1}{N}  \right)  \frac{n}{K} \sum_{k=1}^K \Big(\hat\theta^{(k)} - \hat\theta_{\textup{BAG}}\Big) \Big(\hat\theta^{(k)} - \hat\theta_{\textup{BAG}}\Big) ^\top.
\end{align}
Its properties are given in the following theorem.
See Appendix D for the detailed proof of this theorem.
By Theorem \ref{thm:se}, we can know that $\hat\SE^2( \hat\theta_{\textup{BAG}} )$ should be a reasonable estimator for the asymptotic variance of the bagging estimator $\hat\theta_{\textup{BAG}}$. Thus, we can use it to make further inference, such as constructing confidence interval.
We illustrate the performance of $\hat\SE^2( \hat\theta_{\textup{BAG}} ) $ through extensive simulation experiments.

\begin{theorem}\label{thm:se}
	Assume conditions (C1)--(C4). Further assume that $n\to\infty$, $n/N\to 0$, and $K\to \infty$ as $N\to \infty$. Then we have $ E\Big\{\hat\SE^2(\hat\theta_{\textup{BAG}})\Big\}  =  \big\{(nK)^{-1} + N^{-1} \big\}\Sigma +o \big\{(nK)^{-1} + N^{-1} \big\}$. 
\end{theorem}

\section{NUMERICAL STUDIES}

\subsection{Simulation studies}

To verify the finite performance of the bagging estimator, we conduct a number of simulation studies. 
Assume the whole sample contains a total of $N = 200,000$ observations. 
For each $1\le i\le N$, we generate $(X_i, Y_i)$ under three different models: a linear regression model, a logistic regression model, and a  Poisson regression model.
We set the true parameter as $\theta_0 = (-0.2, -0.1, 0 , 0.1, 0.2)^\top$ for every model.
More detailed model settings are given as follows.

\textbf{Example 1. (Linear Regression)} In this example, we consider $p=5$ covariates $X_i = (X_{i1}, X_{i2}, \dots, X_{i5})^\top$, where each covariate is independently generated from a standard normal distribution. Then, the response $Y_i$ is generated by a linear relationship with the covariates $X_i$ given as 
\begin{align*}
	Y_i = X_i^\top\theta_0 + \varepsilon_i,
\end{align*}
where the noise term $\varepsilon_i$ is independently generated from a standard normal distribution. 

\textbf{Example 2. (Logistic Regression)} In this example, we also consider $p=5$ covariates. As above,  each of them is independently generated from a standard normal distribution. Then, the response $Y_i$ is generated from a Bernoulli distribution with the response probability given as
\begin{align*}
	P(Y_i =1| X_i, \theta_0) = \frac{\exp(X_i^\top \theta_0) }{1+\exp(X_i^\top \theta_0)}.
\end{align*}

\textbf{Example 3. (Poisson Regression)}  In this example, we also consider $p=5$ covariates. They are generated from a multivariate normal distribution $\mN(0, \Sigma_X)$, where $\Sigma_X = (\Sigma_{X,ij})_{1\le i,j\le p} $ with $\Sigma_{X,ij} = 0.5^{|i-j|}$.
Then, the response $Y_i$ is generated from a Poisson distribution given as 
\begin{align*}
	P(Y_i =k| X_i , \theta_0) = \frac{\lambda_i^{k}}{k !} \exp(-\lambda_i), \, \text{where } \lambda_i = \exp(X_i^\top \theta_0).
\end{align*}

%We first consider the logistic regression example. 

Once the full dataset $\{(X_i, Y_i): 1\le i \le N \}$ is simulated, we compute the bagging estimator by using the procedures described above.
For a comprehensive evaluation, various bagging subsample sizes $(n = 500, 750, 1000)$ and numbers of subsamples $(K = 50,100,150,200,250)$ are considered.
Let $\hat \theta_{\textup{BAG}}^{(b)} = \big(\hat \theta_{\textup{BAG},j}^{(b)}  \big)_{j=1}^p$ be the bagging estimator obtained in the $b$-th ($1\le b \le B$) simulation with $B=1,000$. 
We compute the bias as $\Bias_j = |\bar\theta_{\textup{BAG},j} - \theta_{0j} | $, where $\theta_{0j}$ is the $j$-th component of the true value $\theta_0$, and  $ \bar\theta_{\textup{BAG},j}  = B^{-1}  \sum_{b=1}^B \hat\theta_{\textup{BAG},j}^{(b)}$.
The standard error $\SE_j^{(b)}$ can be estimated by using \eqref{eq:se},  that is $\hat\SE_j^{(b)} = \sqrt{\hat S_{jj}^{(b)} }$, where $\hat S_{jj}^{(b)} $ is the $j$-th diagonal element of $\hat\SE^2( \hat \theta_{\textup{BAG}}^{(b)} )$. We report the average $\hat\SE_j = B^{-1} \sum_{b=1}^B \SE_j^{(b)} $.
We also compare  $\hat\SE_j $ with the Monte Carlo standard deviation of $\hat\theta_{\textup{BAG}, j}^{(b)}$. This is calculated by $\SE_j = \sqrt{B^{-1} \sum_{b=1}^B  \big( \hat\theta_{\textup{BAG},j} ^{(b)}  - \bar\theta_{\textup{BAG},j} \big)^2 } $.
Last, we construct a $95\%$ confidence interval for $\theta_j$ as $\CI_j^{(b)} = \big( \hat\theta_{\textup{BAG}, j}^{(b)} - z_{0.975}\hat\SE_j^{(b)},   \hat\theta_{\textup{BAG}, j}^{(b)} + z_{0.975}\hat\SE_j^{(b)}\big)$, where $z_\alpha$ is the $\alpha$-th lower quantile of the standard normal distribution.
Then, the coverage probability is computed as $\ECP_j = B^{-1}\sum_{b=1}^B I \big\{\theta_{0j} \in \CI_j^{(b)}\big\}$, where $I\{ \cdot\}$ is the indicator function. 
The detailed results are given in Tables \ref{tab:linear}--\ref{tab:Poisson} for the linear regression, the logistic regression, and the Poisson regression, respectively. 
% We randomly replicate the experiments a total number of $B=1000$ times.
% To save space, we present the results of logistic regression with $j=2,3$ in Table \ref{tab:logistic}.
% The results of Poisson regression with $j=1,4$ are presented in Table \ref{tab:Poisson}.

In general, the simulation results in Tables \ref{tab:linear}--\ref{tab:Poisson} are very similar. We can draw following conclusions.
First, the SE values steadily decrease as the number of subsamples $K$ increases.
This is expected, because according to Theorem \ref{thm:bias_var}, a larger $K$ leads to a smaller variance.
Note that for the linear regression model, the ordinary least squares estimator is unbiased. Thus, we focus on the estimation bias of the logistic and Poisson regression models.
From Table \ref{tab:logistic} and \ref{tab:Poisson}, we can see that the Bias values does not show a clear tendency to decrease as $K$ increase.
On the other hand, the Bias values decrease as the subsample size $n$ increases.
This is in line with our findings in Theorem \ref{thm:bias_var}, since the averaging operation in bagging does not help to reduce bias.
Last, we find that the estimated $\hat\SE$ values are very close to the Monte Carlo SE values, suggesting that the asymptotic variance formula constructed in \eqref{eq:se} should be correct.
In addition, we can see the empirical coverage probabilities are all around the nominal level 95\%. This confirms the asymptotic normality of the proposed estimator.
It also demonstrates the usefulness of the proposed variance estimator.

\begin{figure}[htbp]
	\centering
	\includegraphics[width = 0.32\textwidth]{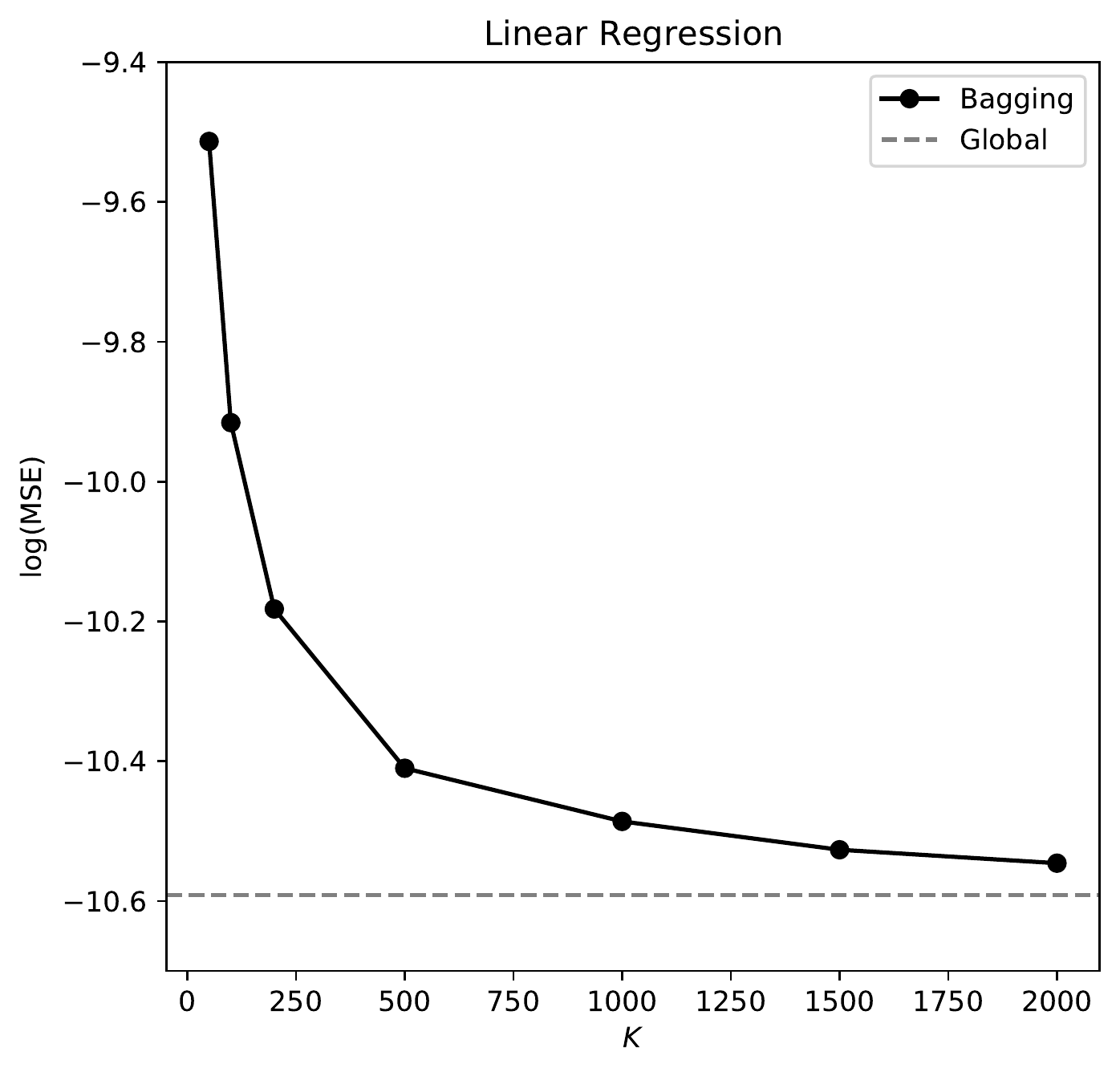}
	\includegraphics[width = 0.32\textwidth]{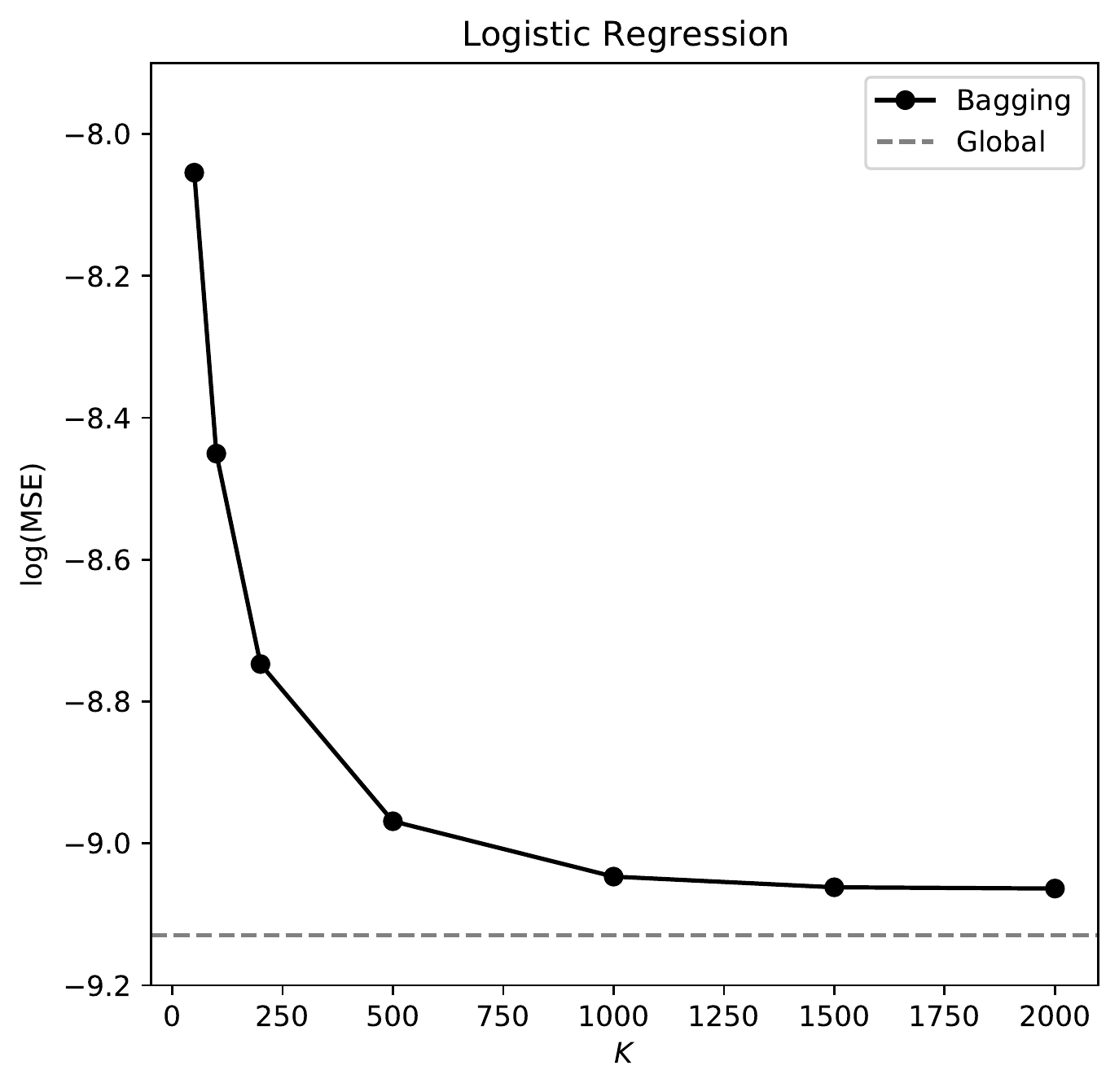}
	\includegraphics[width = 0.32\textwidth]{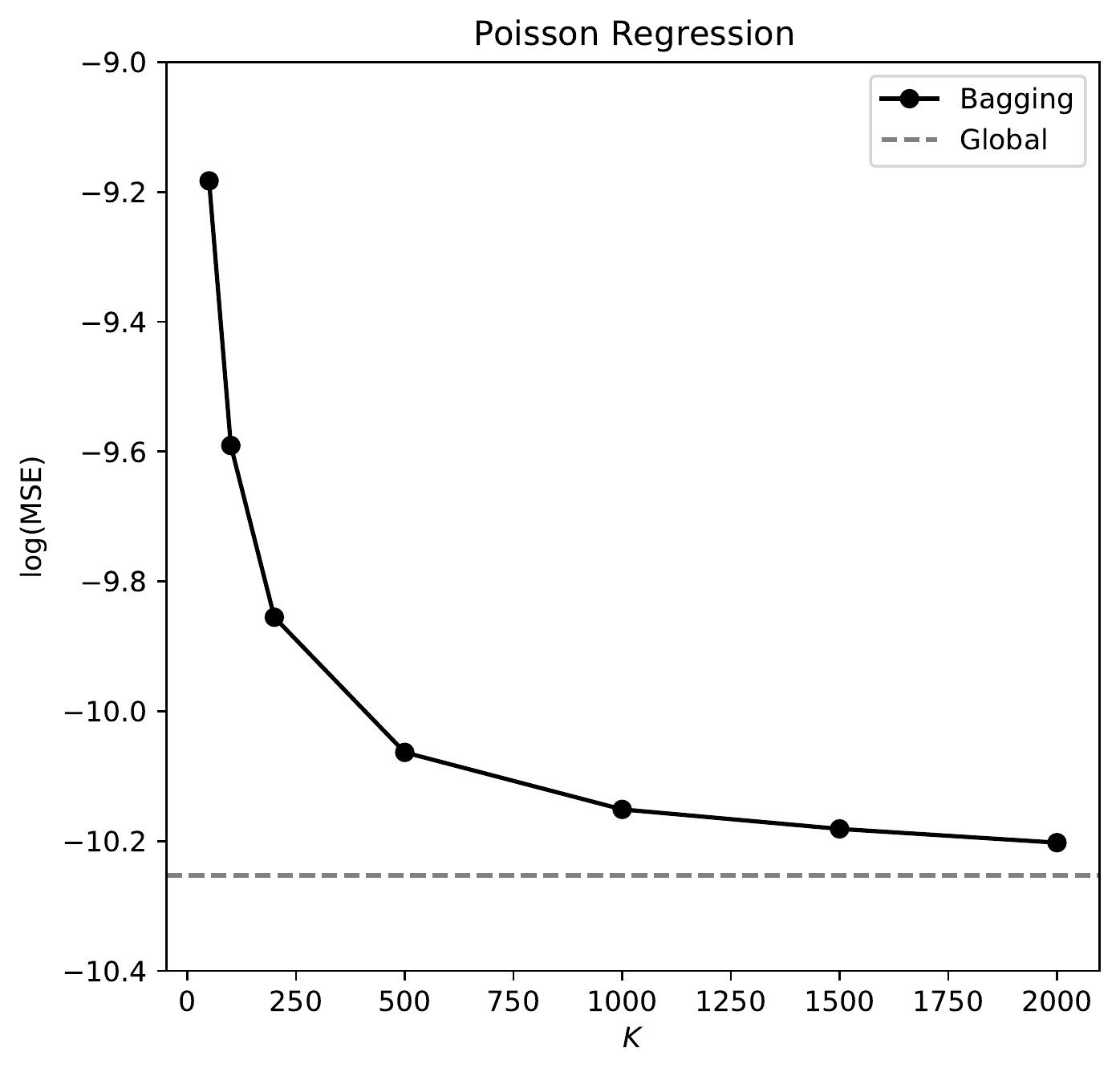}
	\caption{The log-transformed MSE values of the bagging estimator and the global estimator under three different models. We fix subsample size $n$ as $2,000$ and let the number of subsamples $K$ range from $50$ to $2,000$.
		The solid line stands for the bagging estimator and the dashed horizontal line stands for the global estimator.}
	\label{fig:global_and_bag}
\end{figure}

We next compare the performance of the bagging estimator and the global estimator.
Specifically, for computing the bagging estimator, we fix the bagging subsample size as $n=2,000$, and let the number of subsamples $K$ range from $50$ to $2,000$.
Let $\hat\theta^{(b)}$ be one particular estimator (e.g., the global estimator) obtained in the $b$-th ($1\le b\le B$) simulation with $B=500$. 
We compute the mean squared error (MSE) as $\MSE(\hat\theta) =  {B^{-1} \sum_{b=1}^B \| \hat\theta^{(b)} - \theta_0\|^2 } $.
The log-transformed MSE values of different estimators are then plotted in Figure \ref{fig:global_and_bag}.
From Figure \ref{fig:global_and_bag}, we can see that the MSE values of the bagging estimator decrease steadily as $K$ increases.
Furthermore, the bagging estimator performs nearly as well as the global estimator in terms of the MSE value when $K=2,000$.
This is expected, because according to Theorem \ref{thm:anormal}, the bagging estimator should be asymptotically as efficient as the global estimator, provided $n\gg \sqrt{N}$ and $nK \gg N$.

\subsection{A real data example}

For illustration purposes, we apply the bagging method to a real-world dataset, the US Airline Dataset (\url{http://stat-computing.org/dataexpo/2009}).
It consists of the flight arrival and departure details for all commercial flights within the US, from October 1987 to April 2008.
The dataset contains about 123 million records and takes up 12 GB on a hard drive.
After data preprocessing, a total of $N=120,748,239$ records are retained.
Our task is to predict whether a flight is delayed.
We use a binary variable \texttt{Delayed} to denote the delay status of a flight (1 for delayed and 0 for not delayed). 
According to Federal Aviation Administration regulations, we consider a flight delayed if it arrives 15 minutes later than scheduled.
To predict the delay status, four covariates are considered.
The first one is \texttt{Distance}, which describes the distance between the origin and the destination in miles. This numerical variable is standardized to have mean 0 and variance 1.
The second one is the scheduled departure time (\texttt{DepTime}). It is categorized into four time periods: midnight (0:00-7:00), morning (7:00-12:00), afternoon (12:00-18:00), and evening (18:00-0:00).
The third one is \texttt{DayOfWeek}, which is a categorical variable with seven levels (Monday--Sunday).
The last one is \texttt{Month}, which is a categorical variable with twelve levels (January--December).
Next, all the categorical variables are converted to dummy variables.
Note that the first level of each categorical variable is removed to avoid collinearity.
Finally, a total of $p=22$ variables are used as the predictors (including the intercept) in the model.

We then use the logistic regression model to study the delay status of flights.
To compute the bagging estimator, we set the subsample size as $n=[\sqrt{N}\log\log(N)]=32,126$ and draw a total of $K=1,000$ subsamples. 
The detailed estimation results are presented in Table \ref{tab:Airline}.
From Table \ref{tab:Airline}, we have the following interesting findings.
First, the $p$-values corresponding to all covariates are smaller than $0.001$.
This indicates that each of these covariates has a significant influence on the response.
Second, the estimated coefficient for \texttt{Distance} is positive, suggesting that long distance flights are more likely to be delayed.
In terms of the departure time, more delays occur in the afternoon and evening (with coefficients of 1.057 and 0.857, respectively) than in the midnight and morning.
Next, Friday's flights (with coefficient of 0.249) are more often delayed, while Saturday's flights (with coefficient of -0.167) are less likely.
Last, we can see that the estimated coefficients corresponding to February--November are all negative. This implies that January and December are the months with more flight delays.
All those findings are in line with our usual experience.

\section{CONCLUDING REMARKS}

In this paper, we study the bagging estimator for the large-scale $M$-estimation problem.
We first theoretically investigate the asymptotic bias and variance of the bagging estimator.
It reveals that the bagging subsample size should be sufficiently large to make the bias term negligible.
Next, we establish the asymptotic normality for the bagging estimator.
Sufficient conditions are given for the bagging estimator to achieve the optimal asymptotic efficiency.
To support further inference, we also construct a variance estimator.
Finally, we verify the theoretical results by extensive simulation experiments. 
In addition, a real data example is provided to illustrate the usefulness of the bagging estimator.
To conclude the article, we discuss several interesting topics for further study.
First, we assume that the loss function is sufficiently smooth in this paper.
However, some important problems, such as quantile regression and support vector machine, involve non-differentiable loss functions. 
Investigating the theoretical properties of the bagging estimator for these problems is also very meaningful and challenging.
Second, we find that the simple averaging operation in bagging cannot reduce the bias of the final estimator. 
Therefore, it is very interesting to explore some other aggregating strategies.
Last, variable selection methods are widely used in statistical analysis and have received great attention in the past decade.
It is also worthwhile to apply the bagging method to the large-scale variable selection problems.

%\backmatter

%\section*{Acknowledgments}
%This is acknowledgment text~\cite{Elbaum2002}. Provide text here. This is acknowledgment text. Provide text here. This is acknowledgment text. Provide text here. This is acknowledgment text. Provide text here. This is acknowledgment text. Provide text here. This is acknowledgment text. Provide text here. This is acknowledgment text. Provide text here. This is acknowledgment text. Provide text here. This is acknowledgment text. Provide text here. 
%
%\subsection*{Author contributions}
%
%This is an author contribution text. This is an author contribution text. This is an author contribution text. This is an author contribution text. This is an author contribution text. 
%
%\subsection*{Financial disclosure}
%
%None reported.
%
%\subsection*{Conflict of interest}
%
%The authors declare no potential conflict of interests.
%
%
%\section*{Supporting information}
%
%The following supporting information is available as part of the online article:
%
%\noindent
%\textbf{Figure S1.}
%{500{\uns}hPa geopotential anomalies for GC2C calculated against the ERA Interim reanalysis. The period is 1989--2008.}
%
%\noindent
%\textbf{Figure S2.}
%{The SST anomalies for GC2C calculated against the observations (OIsst).}

\appendix

\renewcommand{\theequation}{A.\arabic{equation}}
\setcounter{equation}{0}

\section*{APPENDIX A: PROOF OF LEMMA 1}%\ref{lemma: root n}}

Note that the objective function $\mL_k(\theta)$ is strictly convex in $\theta$.
Consequently, to prove that $\hat\theta^{(k)}$ is $\sqrt{n}$-consistent, it suffices to use the technique in \cite{fan2001variable} to verify that, for any $\epsilon>0$, there exists a finite constant $C>0$ such that, 
\begin{align*}
\liminf_n P\Big\{ \sup_{\| u\| = 1}  \mL_k(\theta_0 + C n^{-1/2} u  ) > \mL_k(\theta_0) \Big\} \ge 1 - \epsilon.
\end{align*}
Write $ \theta_u =\theta_0 + C n^{-1/2} u $, where $C>0$ is a fixed constant and $u\in\mR^p$ is a vector with unit length (i.e., $\| u \| = 1$).
Then by Taylor expansion, we have
\begin{align}
\sup_{\| u\| = 1} \Big\{  \mL_k(\theta_0 + C n^{-1/2} u ) - \mL_k(\theta_0) \Big\} = & n^{-1/2} C u^\top \dot\mL_k(\theta_0)  + (2n)^{-1} C^2 u^\top \ddot\mL_k(\theta_0) u + o_p(1) \nonumber \\
=& C u^\top J_1 + C^2 u^\top J_2 u + o_p(1), \label{eq:quad}
\end{align}
where $J_1 = n^{-1/2} \dot\mL_k(\theta_0)$ and $J_2 = n^{-1}  \ddot\mL_k(\theta_0)$.

We then compute $E(J_1)$ and $\var(J_1)$ as follows. In fact, 
\begin{align*}
E (J_1) = E\{ E(J_1 |\mbS)  \} = \sqrt{n} E \Big\{  E \Big(n^{-1}  \sum_{i\in \mS_k} \ell_i(\theta_0) \Big| \mbS  \Big) \Big\} = \sqrt{n} E \Big\{ N^{-1}  \dot\mL(\theta_0)  \Big\}=0,
\end{align*}
where $\mbS$ denotes the information contained in the whole sample.
In terms of variance, we have $\var(J_1)=$
\begin{align*}
& E\big\{ \var(J_1 | \mbS) \big\} + \var\big\{ E(J_1 | \mbS) \big\} = E\Big\{ n^{-1} \var\Big( \sum_{i\in \mS_k} \ell_i(\theta_0) \Big| \mbS \Big) \Big\} + \var\Big\{ \sqrt{n} N^{-1} \dot \mL(\theta_0)  \Big\}  \\
=& E \left[ N^{-1} \sum_{i = 1}^N  \Big\{ \dot\ell_i(\theta_0) - N^{-1}\dot\mL(\theta)  \Big\} \Big\{ \dot\ell_i(\theta_0) - N^{-1}\dot\mL(\theta)  \Big\}^\top   \right] +  \var\Big\{ \sqrt{n} N^{-1} \dot \mL(\theta_0)  \Big\} \\
=&  \frac{N-1}{N} V + \frac{n}{N} V = O(1).
\end{align*}
Then we should have $J_1 = O_p(1)$. By similar arguments, we can show that $J_2 \to_p H$. 
Recall that $H$ is a positive-definite matrix. Consequently, as long as $C$ is sufficiently large, the quadratic term in \eqref{eq:quad} dominates its linear term. 
This implies that $\sup_{\| u\| = 1} \Big\{  \mL_k(\theta_0 + C n^{-1/2} u ) - \mL_k(\theta_0) \Big\}  > 0$ holds with probability tending to 1 as $n\to \infty$.
This further suggests that, with  probability tending to 1, a local minimizer (i.e., $\hat\theta^{(k)}$) exists, such that $\hat\theta^{(k)} - \theta_0 = O_p(n^{-1/2})$. 
This completes the proof of the lemma.

\renewcommand{\theequation}{B.\arabic{equation}}
\setcounter{equation}{0}

\section*{APPENDIX B: PROOF OF THEOREM 1}%\ref{thm:bias_var}}

We first investigate the term $Q_1$.
It can be computed that $ E(Q_1) = E\big\{ E(Q_1 | \mbS)  \big\}= E\big\{ E(Q_1^{(k)} | \mbS)  \big\} = -H^{-1} E\big\{  N^{-1}\dot\mL(\theta_0) \big\} =0$,
where $\mbS$ denotes the information contained in the whole sample.
We next investigate $\var(Q_1)$. Note that $\var(Q_1) =E\big\{ \var(Q_1 | \mbS) \big\} + \var\big\{ E(Q_1 | \mbS) \big\}$. We can compute that $E \Big\{\var(Q_1 | \mbS) \Big\}=$
\begin{align*}
& K^{-1} H^{-1} E\Big\{ \var\Big(n^{-1} \dot \mL_k(\theta_0) | \mbS\Big)\Big\}H^{-1}
= K^{-1} H^{-1} E\Big\{n^{-1} \var\Big( \dot \ell_i(\theta_0) | \mbS\Big)\Big\}  H^{-1}\\
=& (nK)^{-1} H^{-1} E \bigg[  E \Big\{ \dot\ell_i(\theta_0) \dot\ell_i^\top(\theta_0) \Big| \mbS \Big\} -   E \Big\{ \dot\ell_i(\theta_0) \Big| \mbS \Big\} E \Big\{ \dot\ell_i(\theta_0) \Big| \mbS \Big\} ^\top \bigg] H^{-1}\\
=& (nK)^{-1} H^{-1} \bigg[ E \Big\{ \dot\ell_i(\theta_0) \dot\ell_i^\top(\theta_0)  \Big\}  - E \Big\{\Big( N^{-1}  \mL(\theta_0) \Big) \Big( N^{-1} \mL(\theta_0) \Big)^\top \Big\} \bigg] H^{-1}\\
=&(nK)^{-1} (1-N^{-1}) \Sigma. 
\end{align*}
In addition, we have $\var\big\{ E(Q_1 | \mbS) \big\} =  \var \big\{ H^{-1} N^{-1}\dot\mL(\theta_0) \}  =H^{-1} \var\{N^{-1}\dot\mL(\theta_0) \}H^{-1} = N^{-1}\Sigma$.
Together with these results, we conclude that 
\begin{align*}
\var(Q_1) = (nK)^{-1} (1-N^{-1}) \Sigma + N^{-1}\Sigma = \left( \frac{1}{nK} + \frac{1}{N} \right)\Sigma + o\left( \frac{1}{nK} + \frac{1}{N} \right).
\end{align*}

We next investigate the term $Q_2$.
By Lemmas 18 and 22 and upper bound (B.10)  in \cite{huang2019distributed}, we have $E\| \Delta^{(k)}\|^2  = O(n^{-2})$.  Then we have
\begin{align*}
E\| Q_2\| \le K^{-1} \sum_{k=1}^K E\|Q_2^{(k)} \| = \|H^{-1}\| \cdot E \| \Delta^{(k)}\| = O(n^{-1}).
\end{align*}
Consequently, we should have $Q_2 = O(n^{-1})$.
This completes the proof of the theorem.

\renewcommand{\theequation}{C.\arabic{equation}}
\setcounter{equation}{0}

\section*{APPENDIX C: PROOF OF THEOREM 2}%\ref{thm:anormal}}

Let $\tau = \{ (nK)^{-1} + N^{-1} \}^{1/2} $. By theorem \ref{thm:bias_var} we know that $\tau^{-1} Q_2 \to _p 0$ under the assumed conditions. Thus, it suffices to show that $ \tau^{-1} Q_1 \to_d \mN(0, \Sigma)$. Recall that $Q_1 = H^{-1}  K^{-1} \sum_{k=1}^K \dot\mL_k(\theta_0) $, and $\Sigma = H^{-1} V H^{-1}$. 
Thus,  we turn to show that $U = (\tau K)^{-1} \sum_{k=1}^K \dot\mL_k(\theta_0) \to_d \mN(0, V)$.
This is true if we can show that the characteristic function $f(t)  = E\big\{ \exp(i t^\top U) \} \big\} \to \exp(- t^\top V t / 2) $. 
To this end, we first define $\hat V = \var(\dot\ell_i(\theta_0) | \mbS) = N^{-1} \sum_{i=1}^N \dot\ell_i(\theta_0) \dot\ell_i^\top(\theta_0) - \big\{ N^{-1} \dot\mL(\theta_0) \big\}\big\{ N^{-1} \dot\mL(\theta_0) \big\}^\top$.
Then we can compute that
\begin{align*}
f(t) =& E\Bigg[ \exp\bigg\{ \frac{i t^\top}{\tau K} \sum_{k=1}^K \Big(\dot\mL_k(\theta_0)  - N^{-1}\mL(\theta_0)\Big)   \bigg\}  \exp\bigg\{ \frac{it^\top}{\tau }  N^{-1}\dot\mL(\theta_0)  \bigg\} \Bigg] \\
=&E\Bigg( E\Bigg[\exp\bigg\{ \frac{i t^\top  }{\tau K} \sum_{k=1}^K \Big(\dot\mL_k(\theta_0)  - N^{-1}\mL(\theta_0)\Big)   \bigg\} \Bigg| \mbS \Bigg] \exp\bigg\{ \frac{it^\top}{\tau \sqrt{N} }  N^{-1/2}\dot\mL(\theta_0)  \bigg\} \Bigg) \\
=&E\Bigg( E\Bigg[\exp\bigg\{ \frac{i t^\top \hat V^{1/2} }{\tau nK}  \sum_{k=1}^K\sum_{i\in \mS_k} \hat V^{-1/2} \Big(\dot\ell_i(\theta_0)  - N^{-1}\mL(\theta_0)\Big)   \bigg\} \Bigg| \mbS \Bigg]  \exp\bigg\{ \frac{it^\top}{\tau \sqrt{N} }  N^{-1/2}\dot\mL(\theta_0)  \bigg\} \Bigg)\\
=&E\Bigg[ E\Bigg\{\exp\bigg( \frac{i t^\top \hat V^{1/2} } {\tau \sqrt{nK}}  Z_1   \bigg) \Bigg| \mbS \Bigg\}   \exp\bigg( \frac{it^\top}{\tau \sqrt{N} }  Z_2 \bigg) \Bigg],
\end{align*}
where $Z_1 = (nK)^{-1/2}\sum_{k=1}^K\sum_{i\in \mS_k} \hat V^{-1/2} \Big(\dot\ell_i(\theta_0)  - N^{-1}\mL(\theta_0)\Big) $ and $Z_2 = N^{-1/2}\dot\mL(\theta_0)$.
We investigate the limit of $f(t)$ by considering the following cases.

\textbf{Case 1} ($ nK / N \to \infty$).  In this case, we should have $\tau \sqrt{nK} \to \infty$.
Then we can verify that $\| t^\top \hat V^{1/2}\| / (\tau \sqrt{nK}) \to_p 0$.
Note that $E(Z_1 | \mbS) = 0$ and $\var(Z_1 | \mbS) = I_p$, where $I_p \in \mR^{p \times p}$ is the identity matrix. It follows from the central limit theorem that $Z_1 \to_d \mN(0, I_p)$ conditional on $\mbS$. Then we have $E\Big[\exp\Big\{ i t^\top \hat V^{1/2} Z_1  / \big(\tau \sqrt{nK}\big)     \Big\} \Big| \mbS \Big] \to_p 1$.
Further note that (1) $\tau \sqrt{N} \to_p 1$ and (2) $Z_2 \to_d \mN(0, V)$. It follows that $E \Big[\exp\big\{it^\top Z_2 / (\tau \sqrt{N})\big\} \Big] \to \exp(- t^\top V t /2)$. 
Then by the dominated convergence theorem we conclude that $f(t) \to \exp(-t^\top V t)$.
This completes the proof of Case 1.

\textbf{Case 2} ($ nK / N \to 0$). In this case, we should have $\tau \sqrt{N} \to \infty$. Since $Z_2 \to_d \mN(0, V)$, we have $\exp\big\{it^\top Z_2 / (\tau \sqrt{N})\big\} \to_p 1$. Then by dominated convergence theorem, it remains to show that $E\Big[\exp\Big\{ \big(i t^\top \hat V^{1/2} Z_1 \big) / \big(\tau \sqrt{nK}\big)     \Big\} \Big] \to \exp(-t^\top V t)$. This is true, due to the following reasons: (1) $\tau \sqrt{nK}\to 1$ and $\hat V \to_p V$, implying $\hat V^{1/2}  / \big(\tau \sqrt{nK}\big) \to_p V^{1/2}$. (2) $Z_1 \to_d \mN(0, I_p)$ by the central limit theorem.
This completes the proof of Case 2.

\textbf{Case 3} ($ nK / N \to C$ for some constant $C>0$).  We first decompose $f(t)$ into $f(t) = f_1(t) + f(t) - f_1(t) $ with $f_1(t) = E \Big[ \Gamma_1 \exp\big\{ it^\top V^{1/2} Z_2 / (\tau \sqrt{N}) \big\} \Big]$ and $f(t) - f_1(t) = E \Big[ (\Gamma_1 - \Gamma_2) \exp\big\{ it^\top V^{1/2} Z_2 / (\tau \sqrt{N}) \big\} \Big]$, where $\Gamma_1 = \exp\big\{ -t^\top \hat V t / (2\tau^2 n K) \big\}$ and $\Gamma_2 = E\Big[\exp\big\{ i t^\top \hat V^{1/2} Z_1  / \big(\tau \sqrt{nK}\big)     \big\} \Big| \mbS \Big]$.

Since $\tau \sqrt{nK} \to \sqrt{C+1}$ and $Z_1 \to_d \mN(0, I_p)$, we have $\Gamma_1 - \Gamma_2 \to_p 0$ conditional on $\mbS$. It then follows from the dominated convergence theorem that $f(t)  - f_1(t) \to 0$. Thus, it remains to investigate the limit of $f_1(t)$.
Note that $\hat V\to_p V$ and $\tau^2 nK \to C+1$. Then we have $\Gamma_1 \to_p \exp\big[-t^\top V t / \{2(C+1)\} \big]$. Since $Z_2 \to_d \mN(0, V)$ and $\tau \sqrt{N} \to \sqrt{C^{-1}+1}$, we should have $E \Big[  \exp\big\{ it^\top V^{1/2} Z_2 / (\tau \sqrt{N}) \big\} \Big] \to \exp\big\{  -t^\top V t / (C^{-1} + 1)\big\}$. Together with above results and the dominated convergence theorem, we conclude that $f_1(t) \to  \exp\big[  -t^\top V t / \{2(C+ 1) \}- t^\top V t / \{2(C^{-1} + 1)\} \big] = \exp(-t^\top V t / 2) $. 
This completes the proof of Case 3 and finishes the proof of the theorem.

\renewcommand{\theequation}{D.\arabic{equation}}
\setcounter{equation}{0}

\section*{APPENDIX D: PROOF OF THEOREM 3}%\ref{thm:se}}

Let $c_0 = n\big\{(nK)^{-1} + N^{-1}\big\} $. Recall that $\hat\theta_{\textup{BAG}} = K^{-1} \sum_{k=1}^{K} \hat\theta^{(k)}$. Then we have
\begin{align*}
\hat\SE^2( \hat\theta_{\textup{BAG}} ) =& \frac{c_0}{K} \sum_{k=1}^K \Big(\hat\theta^{(k)} - \hat\theta_{\textup{BAG}}\Big) \Big(\hat\theta^{(k)} - \hat\theta_{\textup{BAG}}\Big) ^\top \\
=&\frac{c_0}{K} \sum_{k=1}^K \Big(\hat\theta^{(k)} - \theta_0 + \theta_0 - \hat\theta_{\textup{BAG}}\Big) \Big(\hat\theta^{(k)}  - \theta_0 + \theta_0 - \hat\theta_{\textup{BAG}} \Big) ^\top \\
=&c_0 \big(A_1 - A_2 \big),		
\end{align*}
where $A_1 =K^{-1} \sum_{k=1}^K \big(\hat\theta^{(k)}  - \theta_0\big)\big(\hat\theta^{(k)}  - \theta_0\big)^\top$ and $ A_2 = \big( \hat\theta_{\textup{BAG}} - \theta_0\big)\big( \hat\theta_{\textup{BAG}} - \theta_0\big)^\top$.
We next compute $E(A_1) $ and $E(A_2)$.

We first compute $E(A_1) $. 
Recall that  $ \hat\theta^{(k)} = \theta_0 - H^{-1} n^{-1}\dot\mL_k(\theta_0)  - H^{-1} \Delta^{(k)}$, where $\Delta^{(k)} = \Big\{ n^{-1}\ddot\mL_k(\theta_0)  - H  \Big\}  (\hat\theta^{(k)} - \theta_0 ) +n^{-1}\Big\{ \int_0^1   \ddot\mL_k\Big((1-t)\theta_0 + t \hat\theta^{(k)} \Big) dt -  \ddot\mL_k(\theta_0)  \Big\}  (\hat\theta^{(k)} - \theta_0)$.
By proof in Appendix B, we can derive that $E\|n^{-1} \dot\mL_k(\theta_0) \|^2 = O(n^{-1})$ and $E\| \Delta^{(k)}\|^2 = O(n^{-2})$.
Then we can compute that $E(A_1 ) =$
\begin{align*}
& E \bigg[ E\Big\{ K^{-1} \sum_{k=1}^K \big(\hat\theta^{(k)}  - \theta_0\big)\big(\hat\theta^{(k)}  - \theta_0\big)^\top \Big| \mbS\Big\}   \bigg] = E\Big\{  \big(\hat\theta^{(k)}  - \theta_0\big)\big(\hat\theta^{(k)}  - \theta_0\big)^\top \Big\}\\
=&  H^{-1} E \Big[  \big\{n^{-1} \dot\mL_k(\theta_0) \big\} \big\{n^{-1} \dot\mL_k(\theta_0) \big\}^\top  -  n^{-1} \dot\mL_k(\theta_0) \Delta^{(k)\top}  -  \Delta^{(k)} \big\{n^{-1}  \dot\mL_k(\theta_0) \big\}^\top  +    \Delta^{(k)}   \Delta^{(k)\top}   \Big] H^{-1}\\
= &n^{-1} \Sigma + O(n^{-3/2}).
\end{align*} 
Similarly, we can compute that $E(A_2) =  E \Big\{\big( \hat\theta_{\textup{BAG}} - \theta_0\big)\big( \hat\theta_{\textup{BAG}} - \theta_0\big)^\top\Big\}=  \big\{(nK)^{-1} + N^{-1} \big\}\Sigma + O(n^{-3/2})$.
Note that $c_0 = o(1)$ under the assumed conditions. Then we have
\begin{align*}
E\Big\{ \hat\SE^2( \hat\theta_{\textup{BAG}} )\Big\} = \left(\frac{1}{nK} + \frac{1}{N} \right) \Sigma + o\left(\frac{1}{nK} + \frac{1}{N} \right).
\end{align*}
This completes the proof.

\bibliographystyle{apalike}
\bibliography{ref.bib}

\clearpage

\begin{sidewaystable}[htbp]
\caption{The simulation results for the bagging estimator under the linear regression model with $B=1,000$ replications. The Bias$_j$ ($\times 10^2$), SE$_j$ ($\times 10^2$), $\hat\SE_j$ ($\times 10^2$) and $\ECP_j$ are reported under different $(n,K)$-specifications.}
\label{tab:linear}
\centering
\begin{tabular}{ r  *{3}{|ccccc}}
\toprule
& \multicolumn{5}{c|}{$n=500$} & \multicolumn{5}{c|}{$n=750$}& \multicolumn{5}{c}{$n=1000$}\\
$K$ & 50 & 100 & 150 & 200 & 250 & 50 & 100 & 150 & 200 & 250 & 50 & 100 & 150 & 200 & 250\\
\hline
Bias$_1$    &    0.017 & 0.010 & 0.006 & 0.005 & 0.003 & 0.016 & 0.006 & 0.003 & 0.003 & 0.005 & 0.009 & 0.005 & 0.002 & 0.008 & 0.009 \\
SE$_1$      &    0.697 & 0.511 & 0.433 & 0.392 & 0.359 & 0.581 & 0.435 & 0.374 & 0.337 & 0.319 & 0.511 & 0.393 & 0.337 & 0.314 & 0.295 \\
$\hat\SE_1$ &    0.663 & 0.500 & 0.429 & 0.388 & 0.361 & 0.557 & 0.427 & 0.372 & 0.341 & 0.321 & 0.495 & 0.385 & 0.340 & 0.316 & 0.300 \\
ECP$_1$     &    0.927 & 0.943 & 0.943 & 0.946 & 0.951 & 0.932 & 0.941 & 0.951 & 0.954 & 0.956 & 0.936 & 0.945 & 0.954 & 0.957 & 0.967 \\
\hline
Bias$_2$    &    0.026 & 0.040 & 0.029 & 0.028 & 0.031 & 0.021 & 0.030 & 0.030 & 0.029 & 0.025 & 0.039 & 0.028 & 0.029 & 0.025 & 0.022 \\
SE$_2$      &    0.681 & 0.515 & 0.439 & 0.391 & 0.355 & 0.579 & 0.439 & 0.366 & 0.338 & 0.313 & 0.515 & 0.390 & 0.338 & 0.310 & 0.294 \\
$\hat\SE_2$ &    0.664 & 0.501 & 0.429 & 0.388 & 0.362 & 0.554 & 0.427 & 0.373 & 0.342 & 0.322 & 0.496 & 0.386 & 0.341 & 0.316 & 0.300 \\
ECP$_2$     &    0.939 & 0.939 & 0.935 & 0.941 & 0.960 & 0.930 & 0.935 & 0.958 & 0.953 & 0.957 & 0.935 & 0.938 & 0.953 & 0.960 & 0.958 \\
\hline
Bias$_3$    &    0.026 & 0.022 & 0.014 & 0.015 & 0.007 & 0.025 & 0.015 & 0.013 & 0.007 & 0.010 & 0.023 & 0.015 & 0.007 & 0.007 & 0.008 \\
SE$_3$      &    0.688 & 0.516 & 0.442 & 0.393 & 0.368 & 0.581 & 0.443 & 0.381 & 0.348 & 0.330 & 0.516 & 0.393 & 0.347 & 0.322 & 0.305 \\
$\hat\SE_3$ &    0.668 & 0.501 & 0.430 & 0.389 & 0.362 & 0.557 & 0.427 & 0.373 & 0.342 & 0.323 & 0.496 & 0.386 & 0.341 & 0.317 & 0.300 \\
ECP$_3$     &    0.939 & 0.945 & 0.937 & 0.939 & 0.954 & 0.929 & 0.936 & 0.949 & 0.951 & 0.950 & 0.940 & 0.940 & 0.949 & 0.950 & 0.953 \\
\hline
Bias$_4$    &    0.006 & 0.011 & 0.006 & 0.011 & 0.015 & 0.003 & 0.006 & 0.011 & 0.017 & 0.016 & 0.013 & 0.010 & 0.017 & 0.014 & 0.011 \\
SE$_4$      &    0.649 & 0.486 & 0.416 & 0.374 & 0.353 & 0.532 & 0.416 & 0.362 & 0.326 & 0.311 & 0.485 & 0.374 & 0.325 & 0.308 & 0.292 \\
$\hat\SE_4$ &    0.663 & 0.499 & 0.428 & 0.388 & 0.361 & 0.554 & 0.426 & 0.372 & 0.341 & 0.322 & 0.493 & 0.385 & 0.340 & 0.316 & 0.300 \\
ECP$_4$     &    0.958 & 0.953 & 0.958 & 0.961 & 0.951 & 0.957 & 0.959 & 0.955 & 0.960 & 0.950 & 0.947 & 0.954 & 0.956 & 0.953 & 0.957 \\
\hline
Bias$_5$    &    0.001 & 0.003 & 0.004 & 0.003 & 0.007 & 0.006 & 0.005 & 0.008 & 0.007 & 0.010 & 0.002 & 0.003 & 0.006 & 0.009 & 0.008 \\
SE$_5$      &    0.679 & 0.518 & 0.451 & 0.401 & 0.366 & 0.575 & 0.451 & 0.383 & 0.343 & 0.325 & 0.518 & 0.401 & 0.343 & 0.322 & 0.305 \\
$\hat\SE_5$ &    0.665 & 0.499 & 0.429 & 0.389 & 0.362 & 0.555 & 0.427 & 0.373 & 0.342 & 0.323 & 0.495 & 0.387 & 0.342 & 0.317 & 0.301 \\
ECP$_5$     &    0.946 & 0.944 & 0.938 & 0.939 & 0.946 & 0.941 & 0.936 & 0.937 & 0.950 & 0.948 & 0.932 & 0.941 & 0.947 & 0.949 & 0.943 \\

\bottomrule 
\end{tabular}
\end{sidewaystable}

\begin{sidewaystable}[htbp]
\caption{The simulation results for the bagging estimator under the logistic regression model with $B=1,000$ replications. The Bias$_j$ ($\times 10^2$), SE$_j$ ($\times 10^2$), $\hat\SE_j$ ($\times 10^2$) and $\ECP_j$ are reported under different $(n,K)$-specifications.}
\label{tab:logistic}
\centering
\begin{tabular}{ r  *{3}{|ccccc}}% >{\centering\arraybackslash}p{0.1\textwidth} *{6}{>{\centering\arraybackslash}p{0.11\textwidth}} }
\toprule
& \multicolumn{5}{c|}{$n=500$} & \multicolumn{5}{c|}{$n=750$}& \multicolumn{5}{c}{$n=1000$}\\
$K$ & 50 & 100 & 150 & 200 & 250 & 50 & 100 & 150 & 200 & 250 & 50 & 100 & 150 & 200 & 250\\
\hline
Bias$_1$    &    0.315 & 0.280 & 0.275 & 0.278 & 0.302 & 0.209 & 0.181 & 0.199 & 0.210 & 0.208 & 0.138 & 0.137 & 0.162 & 0.160 & 0.156 \\
SE$_1$      &    1.391 & 1.056 & 0.921 & 0.806 & 0.755 & 1.172 & 0.915 & 0.775 & 0.720 & 0.681 & 1.046 & 0.799 & 0.717 & 0.667 & 0.637 \\
$\hat\SE_1$ &    1.379 & 1.035 & 0.890 & 0.806 & 0.751 & 1.146 & 0.880 & 0.769 & 0.706 & 0.665 & 1.017 & 0.796 & 0.702 & 0.651 & 0.618 \\
ECP$_1$     &    0.940 & 0.933 & 0.935 & 0.943 & 0.933 & 0.941 & 0.938 & 0.937 & 0.934 & 0.935 & 0.936 & 0.948 & 0.935 & 0.931 & 0.936 \\
\hline
Bias$_2$    &    0.206 & 0.159 & 0.154 & 0.129 & 0.142 & 0.153 & 0.104 & 0.090 & 0.090 & 0.086 & 0.089 & 0.056 & 0.065 & 0.069 & 0.067 \\
SE$_2$      &    1.369 & 1.011 & 0.886 & 0.808 & 0.749 & 1.140 & 0.879 & 0.772 & 0.697 & 0.657 & 1.003 & 0.803 & 0.697 & 0.644 & 0.620 \\
$\hat\SE_2$ &    1.372 & 1.027 & 0.882 & 0.798 & 0.744 & 1.143 & 0.875 & 0.764 & 0.701 & 0.660 & 1.007 & 0.788 & 0.697 & 0.646 & 0.613 \\
ECP$_2$     &    0.944 & 0.951 & 0.933 & 0.942 & 0.938 & 0.942 & 0.931 & 0.942 & 0.946 & 0.945 & 0.945 & 0.941 & 0.950 & 0.944 & 0.951 \\
\hline
Bias$_3$    &    0.006 & 0.021 & 0.011 & 0.004 & 0.019 & 0.020 & 0.013 & 0.024 & 0.018 & 0.018 & 0.023 & 0.004 & 0.018 & 0.013 & 0.016 \\
SE$_3$      &    1.421 & 1.034 & 0.898 & 0.790 & 0.753 & 1.145 & 0.893 & 0.766 & 0.714 & 0.662 & 1.025 & 0.781 & 0.710 & 0.646 & 0.614 \\
$\hat\SE_3$ &    1.357 & 1.022 & 0.878 & 0.795 & 0.741 & 1.132 & 0.872 & 0.761 & 0.698 & 0.657 & 1.005 & 0.784 & 0.694 & 0.643 & 0.611 \\
ECP$_3$     &    0.933 & 0.938 & 0.945 & 0.957 & 0.943 & 0.946 & 0.942 & 0.955 & 0.945 & 0.942 & 0.939 & 0.956 & 0.946 & 0.946 & 0.954 \\
\hline
Bias$_4$    &    0.153 & 0.126 & 0.143 & 0.152 & 0.147 & 0.095 & 0.094 & 0.101 & 0.096 & 0.084 & 0.057 & 0.081 & 0.074 & 0.060 & 0.060 \\
SE$_4$      &    1.405 & 1.068 & 0.917 & 0.825 & 0.754 & 1.190 & 0.911 & 0.780 & 0.710 & 0.671 & 1.057 & 0.815 & 0.708 & 0.651 & 0.607 \\
$\hat\SE_4$ &    1.375 & 1.031 & 0.883 & 0.801 & 0.745 & 1.143 & 0.875 & 0.765 & 0.701 & 0.660 & 1.013 & 0.790 & 0.697 & 0.646 & 0.613 \\
ECP$_4$     &    0.939 & 0.931 & 0.935 & 0.926 & 0.942 & 0.928 & 0.941 & 0.933 & 0.939 & 0.943 & 0.933 & 0.931 & 0.942 & 0.950 & 0.946 \\
\hline
Bias$_5$    &    0.337 & 0.345 & 0.336 & 0.325 & 0.325 & 0.238 & 0.241 & 0.225 & 0.220 & 0.218 & 0.204 & 0.185 & 0.172 & 0.164 & 0.151 \\
SE$_5$      &    1.387 & 1.024 & 0.886 & 0.792 & 0.736 & 1.131 & 0.880 & 0.759 & 0.695 & 0.652 & 1.013 & 0.786 & 0.694 & 0.638 & 0.603 \\
$\hat\SE_5$ &    1.383 & 1.038 & 0.891 & 0.807 & 0.751 & 1.147 & 0.880 & 0.769 & 0.706 & 0.664 & 1.024 & 0.798 & 0.705 & 0.652 & 0.619 \\
ECP$_5$     &    0.936 & 0.944 & 0.922 & 0.937 & 0.920 & 0.942 & 0.932 & 0.947 & 0.940 & 0.938 & 0.944 & 0.949 & 0.947 & 0.944 & 0.946 \\

\bottomrule 
\end{tabular}
\end{sidewaystable}

% \end{landscape}

\begin{sidewaystable}[htbp]
\caption{The simulation results for the bagging estimator under the Poisson regression model with $B=1,000$ replications. The Bias$_j$ ($\times 10^2$), SE$_j$ ($\times 10^2$), $\hat\SE_j$ ($\times 10^2$) and $\ECP_j$ are reported under different $(n,K)$-specifications.}
\label{tab:Poisson}
\centering
\begin{tabular}{ r  *{3}{|ccccc}}% >{\centering\arraybackslash}p{0.1\textwidth} *{6}{>{\centering\arraybackslash}p{0.11\textwidth}} }
\toprule
& \multicolumn{5}{c|}{$n=500$} & \multicolumn{5}{c|}{$n=750$}& \multicolumn{5}{c}{$n=1000$}\\
$K$ & 50 & 100 & 150 & 200 & 250 & 50 & 100 & 150 & 200 & 250 & 50 & 100 & 150 & 200 & 250\\
\hline
Bias$_1$    &    0.143 & 0.135 & 0.120 & 0.115 & 0.109 & 0.120 & 0.086 & 0.072 & 0.076 & 0.082 & 0.082 & 0.061 & 0.058 & 0.062 & 0.060 \\
SE$_1$      &    0.753 & 0.583 & 0.486 & 0.446 & 0.405 & 0.655 & 0.486 & 0.425 & 0.379 & 0.353 & 0.583 & 0.446 & 0.378 & 0.345 & 0.327 \\
$\hat\SE_1$ &    0.736 & 0.551 & 0.473 & 0.428 & 0.399 & 0.614 & 0.470 & 0.410 & 0.377 & 0.355 & 0.547 & 0.427 & 0.377 & 0.349 & 0.332 \\
ECP$_1$     &    0.938 & 0.924 & 0.934 & 0.936 & 0.943 & 0.910 & 0.937 & 0.937 & 0.943 & 0.952 & 0.931 & 0.943 & 0.945 & 0.954 & 0.949 \\
\hline
Bias$_2$    &    0.018 & 0.028 & 0.037 & 0.045 & 0.051 & 0.000 & 0.018 & 0.030 & 0.032 & 0.032 & 0.001 & 0.017 & 0.023 & 0.025 & 0.030 \\
SE$_2$      &    0.834 & 0.630 & 0.547 & 0.491 & 0.456 & 0.727 & 0.548 & 0.472 & 0.429 & 0.404 & 0.632 & 0.492 & 0.428 & 0.399 & 0.375 \\
$\hat\SE_2$ &    0.827 & 0.621 & 0.533 & 0.483 & 0.451 & 0.692 & 0.532 & 0.464 & 0.426 & 0.401 & 0.617 & 0.481 & 0.426 & 0.395 & 0.375 \\
ECP$_2$     &    0.940 & 0.944 & 0.937 & 0.945 & 0.947 & 0.939 & 0.943 & 0.944 & 0.946 & 0.947 & 0.945 & 0.943 & 0.943 & 0.942 & 0.951 \\
\hline
Bias$_3$    &    0.014 & 0.011 & 0.003 & 0.004 & 0.008 & 0.022 & 0.003 & 0.008 & 0.009 & 0.008 & 0.010 & 0.003 & 0.009 & 0.008 & 0.011 \\
SE$_3$      &    0.860 & 0.631 & 0.526 & 0.476 & 0.440 & 0.712 & 0.528 & 0.458 & 0.416 & 0.394 & 0.630 & 0.475 & 0.413 & 0.386 & 0.368 \\
$\hat\SE_3$ &    0.830 & 0.622 & 0.535 & 0.484 & 0.451 & 0.695 & 0.533 & 0.464 & 0.427 & 0.402 & 0.616 & 0.482 & 0.427 & 0.395 & 0.375 \\
ECP$_3$     &    0.934 & 0.946 & 0.955 & 0.948 & 0.943 & 0.939 & 0.949 & 0.943 & 0.950 & 0.957 & 0.938 & 0.944 & 0.953 & 0.952 & 0.949 \\
\hline
Bias$_4$    &    0.045 & 0.026 & 0.044 & 0.051 & 0.050 & 0.002 & 0.028 & 0.029 & 0.031 & 0.030 & 0.000 & 0.025 & 0.021 & 0.025 & 0.026 \\
SE$_4$      &    0.822 & 0.618 & 0.528 & 0.485 & 0.450 & 0.688 & 0.529 & 0.468 & 0.428 & 0.399 & 0.616 & 0.484 & 0.429 & 0.397 & 0.381 \\
$\hat\SE_4$ &    0.830 & 0.623 & 0.535 & 0.484 & 0.451 & 0.698 & 0.534 & 0.465 & 0.427 & 0.402 & 0.618 & 0.482 & 0.426 & 0.395 & 0.375 \\
ECP$_4$     &    0.949 & 0.953 & 0.958 & 0.958 & 0.953 & 0.945 & 0.958 & 0.956 & 0.953 & 0.957 & 0.941 & 0.954 & 0.952 & 0.950 & 0.946 \\
\hline
Bias$_5$    &    0.146 & 0.148 & 0.126 & 0.120 & 0.115 & 0.116 & 0.091 & 0.085 & 0.077 & 0.077 & 0.095 & 0.067 & 0.060 & 0.057 & 0.056 \\
SE$_5$      &    0.761 & 0.551 & 0.473 & 0.422 & 0.391 & 0.622 & 0.474 & 0.409 & 0.369 & 0.356 & 0.553 & 0.423 & 0.370 & 0.350 & 0.328 \\
$\hat\SE_5$ &    0.731 & 0.550 & 0.472 & 0.428 & 0.399 & 0.612 & 0.471 & 0.411 & 0.377 & 0.355 & 0.544 & 0.424 & 0.375 & 0.348 & 0.331 \\
ECP$_5$     &    0.930 & 0.937 & 0.940 & 0.933 & 0.941 & 0.938 & 0.937 & 0.936 & 0.945 & 0.947 & 0.933 & 0.940 & 0.947 & 0.947 & 0.946 \\

\bottomrule 
\end{tabular}
\end{sidewaystable}

\begin{table}[htbp]
\caption{The estimation results of the bagging estimator for the logistic regression model. We report the estimated coefficients (Estimate), standard error ($\hat\SE$) and $p$-values for all covariates.}
\label{tab:Airline}
\centering
\resizebox*{\textwidth}{!}{
\begin{tabular}{llccc | lllccc }
\toprule
Covariates & Levels & Estimate & $\hat\SE (\times 10^2)$ & $p$-Value & Covariates & Levels & Estimate & $\hat\SE (\times 10^2)$ & $p$-Value\\
\hline
\texttt{Intercept} 	&              	& -2.025 & 0.329 & $<0.001$ & \texttt{Month} & February &  -0.064 & 0.243 & $<0.001$ \\
\texttt{Distance}  	&				&  0.122 & 0.046 & $<0.001$ & 				 & March	&  -0.131 & 0.244 & $<0.001$ \\
\texttt{DepTime} 	& Morning      	&  0.396 & 0.261 & $<0.001$ &  				 & April	&  -0.337 & 0.246 & $<0.001$ \\
& Afternoon    	&  0.857 & 0.253 & $<0.001$ &  				 & May		&  -0.341 & 0.249 & $<0.001$ \\
& Evening      	&  1.057 & 0.260 & $<0.001$ &  				 & June		&  -0.025 & 0.240 & $<0.001$ \\
\texttt{DayOfWeek} 	& Tuesday      	& -0.048 & 0.195 & $<0.001$ &  				 & July		&  -0.098 & 0.229 & $<0.001$ \\
& Wednesday    	&  0.046 & 0.194 & $<0.001$ &  				 & August	&  -0.152 & 0.241 & $<0.001$ \\
& Thursday     	&  0.198 & 0.194 & $<0.001$ &  				 & September&  -0.534 & 0.263 & $<0.001$ \\
& Friday       	&  0.249 & 0.186 & $<0.001$ &  				 & October	&  -0.384 & 0.245 & $<0.001$ \\
& Saturday     	& -0.167 & 0.205 & $<0.001$ &  				 & November	&  -0.275 & 0.256 & $<0.001$ \\
& Sunday       	& -0.015 & 0.189 & $<0.001$ &  				 & December	&   0.154 & 0.236 & $<0.001$ \\

\bottomrule

\end{tabular}
}

\end{table}

	%\end{CJK}
\end{document}